\documentclass[12pt]{amsart}
\usepackage[english]{babel}
\usepackage{amsmath}
\usepackage{amsthm}
\usepackage{amssymb}
\usepackage{mathrsfs}
\usepackage{enumerate}
\usepackage[notcite, final, notref]{showkeys}
\usepackage{dsfont}
\usepackage[dvips]{color}
\usepackage{xcolor}
\newtheorem{theorem}{Theorem}[section]

\newtheorem{lemma}[theorem]{Lemma}

\newtheorem{corollary}[theorem]{Corollary}
\newtheorem{definition}[theorem]{Definition}
\theoremstyle{definition}
\newtheorem{remark}[theorem]{Remark}

\def\N{\mathcal{N}}

\def\C{\bf C}
\def\D{\mathcal D}
\def\E{\bf E}
\def\i{\bf i}
\def\R{\mathbb R}
\def\H{\mathbb H}
\def\ux{\underline{x}}
\def\uxi{\underline{\xi}}
\def\uu{\underline{u}}
\def\ut{\underline{t}}
\def\uq{\underline{q}}

\title[]{Dirichlet Spaces In Balls And Half-spaces of $\R^n$}

\author[Y. Yang]{Yan Yang}
\address{Yan Yang, School of Mathematics(Zhuhai)\\
  Sun Yat-sen University\\
 China}
\email{yangyan8@mail.sysu.edu.cn}

\author[T. Qian]{Tao Qian*}
\address{Tao QIAN, Macao Center for Mathematical Sciences\\
Macau University of Science and Technology\\
Macau}
\email{tqian@must.edu.mo}
\thanks{*Corresponding author.\\
Funded by the Science and Technology Development Fund of Macau SAR (grant number 0128/2022/A)}

\begin{document}
\maketitle
\def\e{\bf e}
\begin{abstract}
The present paper studies the Dirichlet spaces in balls and upper-half Euclidean spaces. As main results, we give identical characterizations of the Dirichlet norms in the respective contexts as for the classical 2-D disc case proved by Douglas and Ahlfors.
\end{abstract}

\\

\date{today}

\def\C{\mathbb C}

\section{Introduction}

Denote by ${\mathbb D}$
 the unit disk and $T=\partial {\mathbb D}$ the unit circle in the complex plane $\C.$

Let
\begin{eqnarray*}
h^2({\mathbb D})=\{u(z)\ |  \mbox{ } u(z) \mbox{ is real-valued harmonic in } {\mathbb D} \mbox{ and satisfies }\|u\|_{h^2} <+\infty \}
\end{eqnarray*}
be the real harmonic Hardy space, where $$\|u\|_{h^2}=\sup_{0<r<1}\left(\int_0^{2\pi}|u(re^{{\bf i}t})|^2dt\right)^{\frac{1}{2}}.$$

Define
\begin{eqnarray*}
{\mathcal D}({\mathbb D}):=\{u(z)\ |\  u(z) \mbox { is harmonic in } {\mathbb D}
\mbox{ and satisfies}
\int_{{\mathbb D}}|\nabla u(z)|^2d\sigma<+\infty \}
\end{eqnarray*}
to be the harmonic Dirichlet space in ${\mathbb D},$
where $d\sigma$ is the area measure.

Define the semi-norm of ${\mathcal D}({\mathbb D})$ by
$$\|u\|_{\D({\mathbb D}),*}^2:=\int_{{\mathbb D}}|\nabla  u(z)|^2d\sigma <+\infty $$
and the norm of ${\mathcal D}({\mathbb D})$ by
$$\|u\|_{{\mathbb D}}^2:=\|u\|_{\D({\mathbb D}),*}^2+\|u\|_{h^2}^2.$$
We will call such defined semi-norm as \emph{the gradient form semi-norm}.

If $u\in h^2(\mathbb D)$, Douglas established,  in \cite{Douglas},  the identical relation
\begin{equation}\label{addeqy1}
\int_{\mathbb D} |\nabla u|^2 d\sigma=\frac{1}{2\pi}\int_{T}\int_{T}\frac{|f(z_1)-f(z_2)|^2}{|z_1-z_2|^2}|dz_1| |dz_2|
\end{equation}
that was used as the main technical tool to solve the minimum surface problem of Plateau. He proved the identical relation through a third identical quantity in terms of the Fourier coefficients of the boundary function $f$ (see below Theorem \ref{WZ}), named as \emph{the Fourier series (integral) form} in the sequel. We will call the right-hand-side quantity of (\ref{addeqy1}) as the  \emph{double (singular) integral of the difference-quotient form}, or simply \emph{the double integral form} quantity.

We known that if $u$ is harmonic in $\Omega$, then there exists the corresponding canonical conjugate harmonic function $v$ in $\Omega$ (that is with $v(0)=0$) such that $F=u+{\bf i}v$ is holomorphic in $\Omega$.

Using the Cauchy-Riemann equation, we have
$$F'(z)=\frac{\partial u}{\partial x}+{\bf i}\frac{\partial v}{\partial x}= \frac{\partial u}{\partial x}-{\bf i}\frac{\partial u}{\partial y}.$$
Therefore, we also have
\begin{equation}\label{ca}
\int_{\mathbb D} |\nabla u|^2 dxdy=\int_{\mathbb D} |F'(z)|^2 dxdy\nonumber.
\end{equation}

In \cite{Ahlfors}, using the identity in terms of the canonical holomorphic function
\begin{eqnarray}\label{more}
\int_{\mathbb D} |F'(z)|^2 d\sigma=\lim_{r\rightarrow 1^-}\frac{1}{2}\int_{|z|=r}\overline{F(z)}F'(z)\frac{dz}{\bf i},
\end{eqnarray}
Ahlfors gave another proof of (\ref{addeqy1}). In this paper we will call the left-hand-side of (\ref{more}) as \emph{the holomorphic-derivative form} and the right-hand-side {\it the Ahlfors form quantity}.

In summary, there exist the following equivalent characterizations of $\|u\|_{\D({\mathbb D}),*}$:

\begin{theorem}\label{WZ}$^{\mbox{\scriptsize \cite{WZ, Ahlfors}}}$

If $u\in h^2({\mathbb D})$, $f$ is its boundary function and $F$ is the canonical holomorphic function with the real part $u$, then
\begin{eqnarray*}\label{more1}
\|u\|_{\D({\mathbb D}),*}^2&: =&  \int_{{\mathbb D}} |\nabla u|^2 d\sigma\nonumber\\
&=&\int_{\mathbb D} |F'(z)|^2 d\sigma \nonumber\\
&=& \pi \sum_{k=1}^{+\infty}k(a_k^2+b_k^2)\nonumber\\
&=& \frac{1}{2\pi}\int_{T}\int_{T}\frac{|f(z_1)-f(z_2)|^2}{|z_1-z_2|^2}|dz_1| |dz_2|,\nonumber\\
\end{eqnarray*}
where $\displaystyle a_k=\frac{1}{\pi}\int_0^{2\pi}f(e^{{\bf i}\theta})\cos(k\theta)d\theta$ and $\displaystyle b_k=\frac{1}{\pi}\int_0^{2\pi}f(e^{{\bf i}\theta})\sin(k\theta)d\theta$.
\end{theorem}

We will refer the set relations as \emph{the Douglas-Ahlfors identical relations}. We note that the described quantity comes directly from the real valued harmonic function $u,$ while the one representing it through $F'(z)$ establishes a bridge leading to the complex analysis method used by Ahlfors. See (\ref{more}), also (\ref{AQ}) and (\ref{AH}).\\

\begin{remark}\label{remark1}

Let $u(z)$ be harmonic in ${\mathbb D}$ and $F(z)=u(z)+{\bf i}v(z)$ be the corresponding canonical holomorphic function in ${\mathbb D}.$
Using the Taylor series expansion, we have $F(z)=\displaystyle\sum_{k=0}^{\infty}c_kz^k$. Denote $c_k=a_k-{\bf i}b_k$ and $z=re^{{\bf i}\theta}$,
we obtain $$u(z)=\sum_{k=0}^{\infty}r^k[a_k\cos(k\theta)+b_k\sin(k\theta)], \mbox{ }r<1.$$
If $u(z)\in {\mathcal D}({\mathbb D})$, we have the identical relation between the gradient form and the Fourier series form: $$\|u\|_{\D({\mathbb D}),*}^2=\int_{{\mathbb D}} |\nabla u|^2 d\sigma=\pi \sum_{k=1}^{+\infty}k(a_k^2+b_k^2)<+\infty.$$
Clearly the Dirichlet space ${\mathcal D}({\mathbb D})$ is contained in the harmonic Hardy space $h^2({\mathbb D}),$ then functions belonging to the Dirichlet space all have non-tangential limits a.e. on the boundary of ${\mathbb D}$. Furthermore, there holds
$${\mathcal D}({\mathbb D})=h^2({\mathbb D})\cap {\mathcal H}^{1/2}(T),$$
where ${\mathcal H}^{1/2}(T)$ is the fractional Sobolev space containing functions $f\in L^2(T)$ having the $``1/2"$ order derivative in $L^2(T)$ (\cite{NS}). Precisely,
if $\displaystyle f(e^{{\bf i}\theta})=\sum_{k=0}^{+\infty}[a_k\cos(k\theta)+b_k\sin(k\theta)]\in {\mathcal H}^{1/2}(T),$ then $$\|f\|_{{\mathcal H}^{1/2}(T)}^2=\sum_{k=1}^{+\infty}k(a_k^2+b_k^2)<+\infty.$$
\end{remark}

Studies of Dirichlet space in one complex variable have achieved fruitful results, see, for instance, \cite{ARSW} and \cite{R}, and the references thereby. Recently, the authors of \cite{WZ} studied Dirichlet spaces over bounded Jordan domains of one complex variable. They, in particular, concerned bounded chord-arc domains, for the latter the Douglas-Ahlfors identical relations stand as a characteristic property. The identical relation of the  quantities can be directly interpreted as identical forms of the semi-norm of the Dirichlet spaces on bounded chord-arc domains. The existing studies, however, are restricted to one complex variable. Results analogous to \cite{Douglas} in higher dimensions are obtained in \cite{YQ}. Corresponding to the Dirichlet space in the unit ball, we are interested in whether there are similar results for the upper half plane. The Dirichlet space in the half upper space corresponds to the special case of the general $3$-indices Besov space $\Lambda^{p,q}_\alpha$ for $p=q=2, \alpha =\frac{1}{2}.$ In \cite{S} the space together with its norm are defined through a suitable $p,q$-variation of the double singular integral (see below Theorem \ref{DS}). Among the related results it is known that a function belongs to $\Lambda^{p,q}_\alpha$ if and only if one of some two related quantities is finite. Of the two quantities one is in the Fourier integral form and the other in the gradient form (see also \cite{JP}). The contribution of the present work is to show that for the special case $\Lambda^{2,2}_{\frac{1}{2}}$ the mentioned three quantities are not only equivalent but identical.\\

The present paper studies Dirichlet spaces on balls and upper-half spaces of all dimensions. In proving the identical relations between the quantities of different forms the ones involving the double integral form are most technical.  For the classic unit disc case Douglas proved the identical relations adopting the passage from the gradient form to Fourier series (integral) form, and then to the double integral form; while Ahlfors' adopted the passage from the holomorphic derivative form to the Ahlfors form, via a basic use of the Stokes formula, and then to the double integral form.  Our previous study \cite{YQ} on balls of all dimensions basically goes through Douglas' passage combining with spherical Clifford monogenic expansion and spherical Dirac operator. With the Clifford monogenic expansion treatment the Ahlfors' form quantity comes into play. For unbounded domain contexts, including the upper-half spaces, there are no Fourier form series expansions as Douglas used, and no Ahlfors' form residue computation available. New methods, hence, are to be explored.\\

The writing plan of the paper is as follows: In \S 2 and \S 3 we concentrate in real analysis methods for the harmonic Dirichlet spaces. In \S 2, we define the real- (or harmonic) Dirichlet space ${\mathcal D}(B_{n})$ of the unit balls $B_{n}.$ The gradient energy is used to define the semi-norm of the Dirichlet space. By recalling previously obtained results we give two more identical definitions of, respectively, the Fourier series form and the double integral form \cite{YQ}.  \S 3 is devoted to studying of harmonic Dirichlet spaces of upper-half Euclidean spaces. In such unbounded domain context we establish, using harmonic analysis method, an identical relation between the gradient form and the Fourier integral form. We add a little extra by defining a sequence of $g^\lambda$ functions and show in which $g^0$ and $g^1$ correspond to, respectively, the gradient form quantity and the classical $g(f)$ functions. What is missed out by the real methods in \S 3 is the identical relation with the double singular integral quantity. That will be proved by using the hyper-complex analysis method in \S 4.  \S 4, as the main part of the work,  introduces complex and hyper complex settings to study the complex holomorphic, the quaternionic regular and the Clifford monogenic extensions of the underlying harmonic functions in the respective contexts. With such settings the Ahlfors form quantities come into play.  We prove that for the upper-half complex plane and all the higher dimensional upper-half spaces, and for the quaternionic unit ball, the Ahlfors form quantities are identical with the other three harmonic form quantities. We specially note that the Ahlfors identical relation fails to hold for the unit balls $B_{n+1}(n>1)$ in the non-homogeneous $(n+1)$-dimensional setting.  In the cases the Ahlfors form quantities, nevertheless, are directly equivalent with the underlying Dirichlet space norms, but not with the semi-norms. This, in turn, reflects that the quaternionic space ${\mathbb H}$ is different from imbedding of ${\R}^4$ into the Clifford algebra, as in the former all the hyper-complex $k$-forms for $k>1$ are reduced to those for $k\leq 1.$

\bigskip

\section {The Dirichlet space in the unit ball $B_n$}

\bigskip

This section does not contain essentially new results but links the existing relevant ones with the corresponding Dirichlet space. The section is also to make the whole paper coherent, self-containing and complete.  \\

In the homogeneous setting of $\R^n$ denote by
$B_n=\{\ux \ | \ \mbox{ }|\ux|<1,\ \ux \in \R^n\}$
and
$B_n(r)=\{\ux  \ | \mbox{ }|\ux|<r,\ \ux \in \R^n\}$ the unit ball and the ball with radius $r$ in $\R^n$, respectively.  \\

Denote by $$S^{n-1}=\{{\underline{\xi}} \ | \mbox{ }|\uxi|=1,\ \uxi\in \R^n\}$$ the unit sphere in $\R^n$.

Let the harmonic Hardy space
$$h^2(B_n)=\{u(\ux) \ | \ \mbox{ }u(\ux) \mbox{ is harmonic in } B_n \mbox{ and satisfies }\|u\|_{h^2}<+\infty \},$$
where $$\|u\|_{h^2}=\sup_{0<r<1}\left(\int_{S^{n-1}}|u(r\ux)|^2 d\ux\right)^{\frac{1}{2}}.$$

It is standard knowledge that $u\in h^2(B_n)$ implies that $u$ has non-tangential boundary limits a.e. If denoting by $f(\uxi)$ the boundary limit function, then $f\in L^2(S^{n-1})$, and in the $L^2$-sense $f$ has a spherical Fourier-Laplace series expansion:
$$f(\underline{\xi})=\sum_{k=0}^{+\infty}Y_k(f)(\underline{\xi}),$$
where $$Y_k(f)(\underline{\xi})=c_{n,k}\int_{S^{n-1}}f(\underline{\eta})
P_k^n(\underline{\xi}\cdot\underline{\eta})dS_{\underline{\eta}}$$
stands for the projection of $f$ onto ${\mathcal H}^n_k,$
where $P_k^n$ are the Gegenbauer polynomials and
\begin{eqnarray}\label{eqY4}
c_{n, k}=\frac{1}{\omega_{n-1}}\frac{(n+2k-2)\Gamma(n+k-1)}{(n+k-2)k!\Gamma(n-1)},\nonumber
\end{eqnarray}\label{eq4}
$\displaystyle\omega_{n-1}=\frac{2\pi^{\frac{n}{2}}}{\Gamma(\frac{n}{2})}$ is the surface area of $S^{n-1}$.\\
\def\ux{\underline x}
For $0<r<1$, $\ux=r\uxi,$ there hold the relations
\begin{eqnarray}\label{eq4}
u(\ux)&=&\sum_{k=0}^{+\infty}r^k Y_k(f)(\underline{\xi})\nonumber\\
&=&\sum_{k=0}^{+\infty}r^kc_{n,k}\int_{S^{n-1}}f(\underline{\eta})P_k^n(\underline{\xi}\cdot\underline{\eta})dS_{\underline{\eta}}\nonumber\\
&=&\int_{S^{n-1}}f(\underline{\eta})
\sum_{k=0}^{+\infty}r^kc_{n,k}P_k^n(\underline{\xi}\cdot\underline{\eta})
dS_{\underline{\eta}}\nonumber\\
&=&\frac{1}{\omega_{n-1}}\int_{S^{n-1}}P_r(\underline{\eta}, \underline{\xi})f(\underline{\eta})
dS_{\underline{\eta}},\nonumber \\
&:=&(f\ast P_r)(\underline{\xi})\nonumber,
\end{eqnarray}
where
$$P_r(\underline{\eta}, \underline{\xi})=\frac{1}{\omega_{n-1}}\frac{1-r^2}{|\underline{\eta}-r\underline{\xi}|^n}=
\sum_{k=0}^{+\infty}r^kc_{n,k}P_k^n(\underline{\xi}\cdot\underline{\eta})$$
is the Poisson kernel of $B_n$ at $x=r\underline{\xi}$, which, as a harmonic function by itself, gives rise to the unique harmonic extension $u$  of $f$ into the unit ball $B_n.$ \\

On the other hand, for any $u(\ux)\in h^2(B_n)$, there exist $f\in L^2(S^{n-1})$, such that $u(\ux)=(f\ast P_r)(\underline{\xi})$ and $\|u\|_{h^2}=\|f\|_{L^2(S^{n-1})}$.

\begin{lemma}$^{\mbox{\scriptsize \cite{WL}}}$\label{lemma5}
Let $n\geq 2$, $k\in {\mathcal N}^+$. For any orthonormal base $\{y_1, y_2, \cdots, y_{a^n_k}\}$ of ${\mathcal H}^n_k$ and any $\underline{\xi}, \underline{\eta}\in S^{n-1}$, we have
\begin{eqnarray*}
P^n_k(\underline{\xi}\cdot\underline{\eta})=\frac{1}{c_{n, k}}\sum_{j=1}^{a^n_k}\overline{y_j(\underline{\xi})}y_j(\underline{\eta})
=\frac{1}{c_{n, k}}\sum_{j=1}^{a^n_k}\overline{y_j(\underline{\eta})}y_j(\underline{\xi}).
\end{eqnarray*}
\end{lemma}

\bigskip
In virtue of Lemma \ref{lemma5}, there exists another form of the Fourier-Laplace series:

\begin{eqnarray}
f(\underline{\xi})&=&\sum_{k=0}^{+\infty}c_{n,k}\int_{S^{n-1}}f(\underline{\eta})P_k^n(\underline{\xi}\cdot\underline{\eta})dS_{\underline{\eta}}\nonumber\\
&=&\sum_{k=0}^{+\infty}\int_{S^{n-1}}
f(\underline{\eta})\sum_{j=1}^{a^n_k}\overline{y_j(\underline{\eta})}y_j(\underline{\xi})dS_{\underline{\eta}}\nonumber\\
&=&\sum_{k=0}^{+\infty}\sum_{j=1}^{a^n_k} b_j y_j(\underline{\xi})\nonumber,
\end{eqnarray}
where
\begin{equation}\label{eqY14}
b_j=\int_{S^{n-1}}f(\underline{\eta})\overline{y_j(\underline{\eta})}dS_{\underline{\eta}}.
\end{equation}

Clearly, if $f\in L^2(S^{n-1})$, then $\displaystyle\|f\|_2^2=\sum_{k=0}^{+\infty}\sum_{j=1}^{a^n_k} |b_j|^2$.

The Dirichlet space of $B_n$ is defined as

\begin{definition}
Define the Dirichlet space in $B_n$
\begin{eqnarray*}
{\mathcal D}(B_n):=\{u(\ux)\ |  \mbox{ }u(\ux) \mbox{ is harmonic in } B_n
\mbox{ and satisfies }
\int_{B_n}|\nabla u(\ux)|^2d\ux<+\infty \}.
\end{eqnarray*}

Define the semi-norm in ${\mathcal D}(B_n)$
$$\|u\|_{\D(B_n),*}^2:=\int_{B_n}|\nabla u(\ux)|^2d\ux <+\infty $$
and its norm by
$$\|u\|_{\D(B_n)}^2:=\|u\|_{\D(B_n),*}^2+\|u\|_{h^2}^2.$$
\end{definition}

\bigskip
In \cite{YQ} we obtained several equivalent characterizations of $\|u\|_{\D(B_n),*}:$

\begin{theorem}\label{yyth2.1}

Let $u(\ux)\in h^2(B_n)$, $f(\underline{\eta})$ be its boundary function. Then we have
\begin{eqnarray}\label{eqy3}
\|u\|_{\D(B_n),*}^2&:=&\int_{B_n}|\nabla u(\ux)|^2d\ux\nonumber\\
&=&\sum_{k=1}^{+\infty}k\sum_{j=1}^{a^n_k}|b_j|^2\nonumber\\
&=& \frac{1}{\omega_{n-1}}\int_{S^{n-1}}\int_{S^{n-1}}
\frac{|f(\underline{\eta_1})-f(\underline{\eta_2})|^2}{|\underline{\eta_1}-\underline{\eta_2}|^n}
 dS_{\underline{\eta_1}}dS_{\underline{\eta_2}}\nonumber,
\end{eqnarray}
where $b_j$'s are the coefficients of the Fourier-Laplace series expansion of $f$ which are given by the formula (\ref{eqY14}).

\end{theorem}

\begin{remark}

The definition of the Dirichlet space gives rise to the relation
$${\mathcal D}(B_n)=h^2(B_n)\cap{\mathcal H}^{1/2}(S^{n-1}),$$
where ${\mathcal H}^{1/2}(S^{n-1})$ is the fractional Sobolev space containing functions $f\in L^2(S^{n-1})$ having $``1/2"$ derivative in $L^2(S^{n-1})$. More precisely, that is
$$\|f\|_{{\mathcal H}^{1/2}(S^{n-1})}^2=\sum_{k=1}^{+\infty}k\sum_{j=1}^{a^n_k}|b_j|^2<\infty.$$
\end{remark}

\section {Dirichlet spaces in upper-half spaces: The real harmonic analysis methods }

Denote by $\R_+^{n+1}=\{(\ux, y)| \mbox{ }\ux\in \R^n, \ y>0 \}$ the upper-half space in the $(n+1)$-dimensional Euclidean space. \\

Denote the harmonic Hardy $h^2$ space of the half Euclidean space by $$h^2(\R_+^{n+1})=\{u(\ux, y)\ |\ \mbox{ }u(\ux, y) \mbox{ is harmonic in } \R_+^{n+1} \mbox{ and satisfies }\|u\|_{h^2(\R_+^{n+1})}<+\infty \},$$
where $$\|u\|_{h_2(\R_+^{n+1})}=\sup_{0<y<+\infty}\left(\int_{\R^n}|u(\ux, y)|^2 d\ux\right)^{\frac{1}{2}}.$$

It is well-known that a harmonic function $u\in h^2(\R_+^{n+1})$ if and only if there exists $f(\ux)\in L^2{(\R^n)}$ such that
\begin{equation}\label{eqy1}
u(\ux, y)=\int_{\R^n}P_y(\ux-\ut)f(\ut)d\ut,
\end{equation}
where $P_y(\ux-\ut)$ is the Poisson kernel. Moreover, $u$ is harmonic in $\R_+^{n+1},$
$$\lim_{y\rightarrow 0^+}u(\ux, y)=f(\ux) \mbox{ }a.e.,$$
and the Hardy space norm $\|u\|_{h^2}=\|f\|_{L^2}.$ The Poisson kernel in $\R_+^{n+1}$
has the expressions $$P_y(\ux-\ut)=\frac{2}{\omega_n}\frac{y}{(y^2+|\ux-\ut|^2)^{\frac{n+1}{2}}}=\int_{\R^n}e^{-2\pi {\bf i}\uu\cdot (\ux-\ut)}e^{-2\pi|\uu|y}d\uu$$
with $$\frac{2}{\omega_n}=\frac{\Gamma(\frac{n+1}{2})}{\pi^{\frac{n+1}{2}}}.$$

Using the Plancherel Theorem to the right-hand-side of (\ref{eqy1}), we have
\begin{equation}\label{eqy2}
u(\ux, y)=\int_{\R^n}\hat{f}(\ut)e^{-2\pi {\bf i} \ut\cdot \ux}e^{-2\pi|\ut|y}d\ut,\mbox{ }y>0,
\end{equation}
where $\hat{f}$ is the Fourier transform of $f$ defined by
$$\hat{f}(\ut)=\int_{\R^n}e^{2\pi {\bf i} \ux \cdot \ut}f(\ux)d\ux.$$

Next, we define the Dirichlet space in $\R_+^{n+1}$.

\begin{definition}
Define the Dirichlet space
\begin{eqnarray*}
{\mathcal D}(\R_+^{n+1}):=\{u(\ux, y)| \mbox{ }u(\ux, y)\in h^2(\R_+^{n+1})
\mbox{ and satisfies } \int_{\R_+^{n+1}}|\nabla u(\ux, y)|^2d\ux dy<+\infty \}.
\end{eqnarray*}

Define the semi-norm in ${\mathcal D}(\R_+^{n+1})$
$$\|u\|_{\D(\R_+^{n+1}),*}^2:=\int_{\R_+^{n+1}}|\nabla u(\ux, y)|^2d\ux dy<+\infty $$
and the norm
$$\|u\|_{\D(\R_+^{n+1})}^2:=\|u\|_{\D(\R_+^{n+1}),*}^2+\|u\|_{h^2}^2.$$
\end{definition}

\bigskip
Next, we prove the identical relation between $\|u\|_{\D(\R_+^{n+1}),*}$ and the Fourier integral quantity.

\begin{theorem}\label{yth2.1}

Let $u(\ux, y)\in {\mathcal D}(\R_+^{n+1})$, then we have
\begin{eqnarray}\label{eqy3}
\|u\|_{\D(\R_+^{n+1}),*}^2&:=&\int_{\R_+^{n+1}}|\nabla u(\ux,y)|^2d\ux dy\nonumber\\
&=&2\pi\int_{\R^{n}}|\ut||\hat{f}(\ut)|^2 d\ut\nonumber.
\end{eqnarray}
\end{theorem}

{\bf Proof: } By the Fubini Theorem, we have
\begin{eqnarray*}
\int_{\R_+^{n+1}}|\nabla u(\ux, y)|^2d\ux dy=\int_{0}^{+\infty}dy\int_{\R^n}|\nabla u(\ux,y)|^2d\ux.
\end{eqnarray*}
From (\ref{eqy2}), we have
$$\frac{\partial u}{\partial y}=\int_{\R^n}-2\pi|\ut|\hat{f}(\ut)e^{-2\pi {\bf i}\ut\cdot \ux}e^{-2\pi|\ut|y}d\ut$$
and
$$\frac{\partial u}{\partial x_j}=\int_{\R^n}-2\pi {\bf i} t_j \hat{f}(\ut)e^{-2\pi {\bf i}\ut\cdot \ux}e^{-2\pi|\ut|y}d\ut.$$

Therefore, we have
\begin{equation}
\int_{\R^n}|\nabla u(\ux, y)|^2d\ux=\int_{\R^n}8\pi^2|\ut|^2|\hat{f}(\ut)|^2e^{-4\pi|\ut|y}d\ut, \mbox{ }y>0.\nonumber
\end{equation}

Then,
\begin{eqnarray*}
\int_{\R_+^{n+1}}|\nabla u(\ux, y)|^2d\ux dy&=&\int_{0}^{+\infty}dy\int_{\R^n}|\nabla u(\ux,y)|^2d\ux\\
&=&\int_{0}^{+\infty}dy\int_{\R^n}8\pi^2|\ut|^2|\hat{f}(\ut)|^2e^{-4\pi|\ut|y}d\ut\\
&=&\int_{\R^n}8\pi^2|\ut|^2|\hat{f}(\ut)|^2d\ut\int_{0}^{+\infty}e^{-4\pi|\ut|y}dy\\
&=&2\pi\int_{\R^n}|\ut||\hat{f}(\ut)|^2d\ut.
\end{eqnarray*}

This completes the proof.

\begin{remark}
From Theorem \ref{yth2.1}, we have
$${\mathcal D}(\R_+^{n+1})=h^2(\R_+^{n+1})\cap{\mathcal H}^{1/2}(\R^n).$$
where ${\mathcal H}^{1/2}(\R^n)$ is the fractional Sobolev space containing functions $f\in L^2(\R^n)$ having $``1/2"$ derivative in $L^2(\R^n)$. More precisely, that is
$$\|f\|_{{\mathcal H}^{1/2}(\R^n)}^2=\int_{\R^n}|\ut||\hat{f}(\ut)|^2 d\ut.$$
\end{remark}

\bigskip
In the following, we will present the relationship between the Dirichlet space and the $g$-function.

First, let recall the definition of the $g$-function.
\begin{definition}$^{\mbox{\scriptsize \cite{S}}}$
Let $f\in L^2(\R^n)$ be the boundary function of harmonic function $u$ in the upper half space $\R_+^{n+1}$. Define the $g$-function
$$g(f)(\ux):=(\int_0^{+\infty}|\nabla u(\ux, y)|^2 y dy)^{\frac{1}{2}}.$$
\end{definition}
It is known that (\cite{S}) $$\|g(f)\|_2=\frac{\sqrt{2}}{2}\|f\|_2.$$

We introduce a series $g^{\lambda}$-functions ($0\leq\lambda\leq 1$),
where $g^1(f)(\ux)=g(f)(\ux),$ the classical $g$ function of $f$, and
\[ g^0(f)(\ux)=\left(\int_0^{+\infty}|\nabla u(\ux, y)|^2 dy\right)^{\frac{1}{2}}.\]

\begin{definition}
Let $f\in L^2(\R^n)$. For $0\leq\lambda\leq 1$, define the $g^{\lambda}$-function
$$g^\lambda(f)(\ux):=(\int_0^{+\infty}|\nabla u(\ux, y)|^2 y^{\lambda} dy)^{\frac{1}{2}}.$$
\end{definition}

Then we have
\begin{theorem}\label{yyth2.5}
For $0\leq\lambda\leq 1$, if $f\in L^2(\R^n)\cap{\mathcal H}^{(1-\lambda)/2}(\R^n)$, then we have
$$\|g^{\lambda}(f)\|_2=\sqrt{\Gamma(\lambda+1)2^{1-2\lambda}\pi^{1-\lambda}}\|f\|_{{\mathcal H}^{(1-\lambda)/2}(\R^n)},$$
\end{theorem}
where $$\|f\|_{{\mathcal H}^{(1-\lambda)/2}(\R^n)}^2=\int_{\R^n}|\ut|^{(1-\lambda})|\hat{f}(\ut)|^2 d\ut.$$

{\bf Proof: }
\begin{eqnarray*}
\|g^{\lambda}(f)\|_2^2&=&\int_{\R^n}(\int_0^{+\infty}|\nabla u(\ux, y)|^2 y^{\lambda}dy)d\ux\\
&=&\int_{0}^{+\infty} y^{\lambda}dy \int_{\R^n}|\nabla u(\ux,y)|^2d\ux\\
&=&\int_{0}^{+\infty}y^{\lambda} dy\int_{\R^n}8\pi^2|\ut|^2|\hat{f}(\ut)|^2e^{-4\pi|\ut|y}d\ut\\
&=&\int_{\R^n}8\pi^2|\ut|^2|\hat{f}(\ut)|^2d\ut\int_{0}^{+\infty}e^{-4\pi|\ut|y}y^{\lambda}dy\\
&=&\Gamma(\lambda+1)2^{1-2\lambda}\pi^{1-\lambda}\int_{\R^n}|\ut|^{1-\lambda}|\hat{f}(\ut)|^2d\ut.
\end{eqnarray*}
This completes the proof.

From Theorem \ref{yyth2.5}, we obtain
\begin{corollary}
When $\lambda=1$, $g^1(f)(\ux)=g(f)(\ux)$, which is the $g$-function, then
$$\|g^{1}(f)\|_2^2=\|g(f)\|_2^2=\frac{1}{2}\|\hat{f}\|_2^2=\frac{1}{2}\|f\|_2^2.$$

When $\lambda=0$, we have $g^0(f)(\ux)=(\int_0^{+\infty}|\nabla u(\ux, y)|^2 dy)^{\frac{1}{2}}$, then
$$\|g^0(f)\|_2^2=2\pi\|U\|_{\D(\R_+^{n+1}),*}^2=2\pi\|\sqrt{|\ut|}\hat{f}(\ut)\|_2^2=2\pi\|f\|_{{\mathcal H}^{1/2}(\R^n)}^2.$$
\end{corollary}

\section{The complex and hyper-complex settings for the Ahlfors form quantity }
\def\H{\mathbb{H}}
\def\HC{{H}^2(\mathbb{C}^+)}
\def\R{\mathbb{R}}
\def\N{\mathbb{N}}
\def\C{\mathbb{C}}
\def\e{\bf e}
\def\E{\bf E}
\def\i{\bf i}
\def\j{\bf j}
\def\k{\bf k}

\S 4 will be divided into four subsections, establishing the theories on the upper-half complex plane, on the upper-half Euclidean spaces ${\R}^{n+1}_+,$ and on the unit ball and the upper-half Hamilton quaternionic space $\H.$  In the last subsection we give a revision of the balls cases in ${\R}^{n+1}, n>1.$

\subsection{The Dirichlet space in the upper-half complex plane}

\

\bigskip

It is standard knowledge that there exist Moebius transformations that conformally map the unit disk ${\mathbb D}$ to the upper half complex plane ${\mathbb C}_+$.

Let $$z=L(\zeta)=\frac{\zeta-{\bf i}}{\zeta+{\bf i}},$$ that conformally maps the upper half plane ${\mathbb C}_+$ to the unit disk ${\mathbb D},$ and maps their boundaries in one-to-one and onto manner keeping the orientation.

Denote $L(x_1)=z_1$ and $L(x_2)=z_2$.  By direct computation we have
$$|L(x_1)-L(x_2)|=\left|\frac{2{\bf i}(x_1-x_2)}{(x_1+{\bf i})(x_2+{\bf i})}\right|$$
and $$|dL(x)|=\frac{2}{|x+{\bf i}|^2}dx.$$

Let $u\in h^2({\mathbb D}),$ $f$ is its boundary function, then we obtain
\begin{eqnarray*}
&&\frac{1}{2\pi}\int_{T}\int_{T}\frac{|f(z_1)-f(z_2)|^2}{|z_1-z_2|^2}|dz_1| |dz_2|\\
&=&\frac{1}{2\pi}\int_{-\infty}^{+\infty}\int_{-\infty}^{+\infty}\frac{|f(L(x_1))-f(L(x_2))|^2}{|L(x_1)-L(x_2)|^2}|dL(x_1)||dL(x_2)|\\
&=&\frac{1}{2\pi}\int_{-\infty}^{+\infty}\int_{-\infty}^{+\infty}\frac{|\tilde{f}(x_1)-
\tilde{f}(x_2)|^2}{|x_1-x_2|^2}dx_1dx_2,
\end{eqnarray*}
where $\tilde{f}=f(L)=f\circ L$.

Since $L(\zeta)$ is holomorphic, through simple computation, we also have
\begin{eqnarray*}
\int_{{\mathbb D}} |\nabla_{(x,y)} u(z)|^2 dxdy&=&\int_{{\mathbb C}_+} |\nabla_{(s, t)} u(L(\zeta))|^2 dsdt\\
&=&\int_{{\mathbb C}_+} |\nabla_{(s, t)} \tilde{u}(\zeta)|^2 dsdt,
\end{eqnarray*}
where $ \tilde{u}(\zeta)=u(L(\zeta))=u\circ L(\zeta)$ is harmonic in ${\mathbb C}_+$ when $u(z)$ is harmonic in ${\mathbb D}$.

Define
\begin{eqnarray*}
{\mathcal D}({\mathbb C}_+):=\{u(z)\ |\ \mbox{ }u(z) \in h^2({\mathbb C}_+)
\mbox{ and satisfies }
\int_{{\mathbb C}_+}|\nabla u(z)|^2d\sigma<+\infty \}
\end{eqnarray*}
the Dirichlet space in ${\mathbb C}_+$.

Define the semi-norm in ${\mathcal D}({\mathbb C}_+)$
$$\|u\|_{\D({\mathbb C}_+),*}^2:=\int_{{\mathbb C}_+}|\nabla u(z)|^2d\sigma <+\infty $$
and the norm
$$\|u\|_{{\mathbb C}_+}^2:=\|u\|_{\D({\mathbb C}_+),*}^2+\|u\|_{h^2}^2.$$

Summarizing what just obtained and recalling Corollary \ref{last}, we have
\begin{theorem}\label{eqyy1}
 Let $u\in h^2({\mathbb C}_+),$ $f$ is its boundary function, and $F$ is the canonical holomorphic extension into ${\mathbb C}_+$ with real part $u.$ Then
\begin{eqnarray}\label{more1}
\|u\|_{\D({\mathbb C}_+),*}^2&:=&\int_{{\mathbb C}_+} |\nabla u|^2 d\sigma\nonumber\\
&=&2\pi\int_{-\infty}^{+\infty}|t||\hat{f}(t)|^2 dt\nonumber \\
&=&\int_{{\mathbb C}_+} |F'(z)|^2 d\sigma\nonumber\\
&=& \frac{1}{2\pi}\int_{-\infty}^{+\infty}
\int_{-\infty}^{\infty}\frac{|f(x_1)-f(x_2)|^2}{|x_1-x_2|^2}dx_1dx_2,\nonumber
\end{eqnarray}
\end{theorem}
where $\hat{f}$ is the Fourier transform of $f$ defined by $$\hat{f}(t)=\int_{-\infty}^{\infty}e^{2\pi {\bf i} xt}f(x)dx.$$

\begin{remark}
It is known that
Moebius transformations do not, in general, map $L^2$ functions on the unit circle to $L^2$ functions on the line. This suggests that the Dirichlet space $\D({\mathbb C}_+)$ we wish to define should first be a subset of the harmonic Hardy space $h^2({\mathbb C}_+).$ On the other hand, we comment that in the unit disk case this $h^2$ requirement is not so essential, for in the Fourier series case convergence of the weighted Fourier series $\sum^\infty_{k=1}k|b_k|^2$ implies convergence of the original Fourier series of $f,$ or $\|f\|_2^2<\infty.$ The analogous Fourier integral expression does not have this property (see the integral expression involving $|t||\hat{f}(t)|^2$ in Theorem \ref{eqyy1}). In summary, conformal mappings may only preserve the semi-norms, but not the norm. \\

We therefore have
$${\mathcal D}({\mathbb C}_+)=h^2({\mathbb C}_+)\cap{\mathcal H}^{1/2}({\mathbb R}).$$
where ${\mathcal H}^{1/2}({\mathbb R})$ is the fractional Sobolev space containing functions $f\in L^2({\mathbb R})$ having $``1/2"$ derivative in $L^2({\mathbb R})$. More precisely, that is
$$\|f\|_{{\mathcal H}^{1/2}({\mathbb R})}^2=\int_{-\infty}^{+\infty}|t||\hat{f}(t)|^2 dt<\infty.$$
\end{remark}

\begin{remark}
We are to establish the theory in the upper-half spaces ${\R}^{n+1}_+, n\ge 1.$ It is noted that the conformal mapping method to transform a higher-dimensional unit balls into an upper-half space is unavailable. There are at least two obstacles. First, in the higher dimensions, the composition of a Moebius transform and a harmonic function is no longer a harmonic function. Secondly. Moebius transforms in higher dimensions are, by themselves, no longer monogenic (Clifford holomorphic) (\cite{Ahl}). We, therefore, seek for a direct approach to proving analogous results for higher dimensions.
\end{remark}

\def\bi{\bf i}
\subsection{The Dirichlet space in ${\R}^{n+1}_+$}

\
\bigskip

We will give a short revision on Clifford analysis. The basic references include \cite{BDS, GM, DSS}.  Clifford analysis is an extension of analysis of one complex variable to high dimensions. Essentially speaking, Clifford analysis is to add one new real variable $x_0$ to the existing $n$ real variables $x_1,\cdots,x_n,$ so to create a Cauchy structure. In contrast, analysis of several complex variables corresponds to adding $n$
complex parts $y_j{\bi}$ to each existing real variable $x_j, j=1,\cdots,n,$ to form $n$ complex variables
$z_j=x_j+y_j{\bi},$ and so to create a Cauchy structure. The latter is, in fact, real $2n$-dimension. We will study functions defined in one and multi-dimensional Euclidean spaces taking values in real (or complex) Clifford algebra. Similar to complex analysis, the Cauchy integral formula, Taylor and Laurent series expansions etc., the most basic objects in one complex variable, are all available in Clifford analysis. Just as the relationship between holomorphic functions of one complex variable and 2-D harmonic functions, Clifford analysis has close connections with multi-variable harmonic analysis. As examples, a conjugate harmonic system (\cite{SW}) is just the components of a Clifford monogenic function (holomorphic function in Clifford analysis), and Hilbert transform of a function in Euclidean space is the ${\bf e}_j$-multiple-sum of the corresponding Riesz transforms (\cite{DMQ}). Clifford algebraic structure has been well adopted in contemporary harmonic analysis (\cite{GM}, \cite{QY}), and helps to solve deep real-analysis problems (\cite{LMcQ,LMcS,Q}). We prove the identical relations between the double integral form and the other forms in the Douglas-Ahlfors format through the Ahlfors form quantity which gives one more example of the crucial role of Clifford analysis. \\

Let ${\e}_1, {\e}_2, \cdots , {\e}_{n}$ be the basic elements satisfying ${\e}_j{\e}_k+{\e}_k{\e}_j=-2\delta_{jk}$, where $\delta_{jk}=1$ if $j=k$, and $\delta_{jk}=0$ otherwise, $j, k=1, 2, \cdots , n$. Let
$${\R}^n=\{\ux\ |\ \ux=x_1{\e}_1+\cdots +x_{n}{\e}_{n}:x_k \in {\R}, k=1, 2, \cdots , {n}\}$$
be the n-dimensional homogeneous Euclidean space, and

$${\R}^{n+1}=\{(\ux, x_0):=x_0+\ux\ |\ \ux=x_1{\e}_1+\cdots +x_{n}{\e}_{n}:x_k \in {\R}, k=0, 1, \cdots , {n}\}$$
the non-homogeneous Euclidean $(n+1)$-dimensional space we concern.

Let ${\mathcal{CL}}_{0, n}$ denote the real (or complex)
Clifford algebra generated by ${\e}_1, {\e}_2, \cdots , {\e}_{n}$. The linear basis for the
Clifford algebra is given by ${\e}_A$, where $A$ runs over all the
ordered subsets of $\{0, 1, \cdots, n\}$, namely,
 $$A=\{1\leq i_1<i_2<\cdots <i_l\leq n\}, \mbox{ }1\leq l\leq n, $$
with ${\e}_0={\e}_{\emptyset}=1.$

A general element $x$ of ${\mathcal{CL}}_{0, n}$ can be represented in the form $x=\sum_{k=0}^{n-1}[x]_k$, where
$[x]_k=\sum\limits_{A} x_A {\e}_A, {\e}_A={\e}_{i_1}{\e}_{i_2}\cdots {\e}_{i_k}, 1\leq i_1< i_2<\cdots< i_k\leq {n}. \mbox{ }x_A\in \R \ (\mbox{or}\ {\mathbb C})$ and $[x]_0\in \R \ (\mbox{or}\ {\mathbb C}).$

Let $x=\sum_{k=0}^{n}[x]_k \in {\mathcal{CL}}_{0, n}$, then $x$ consists of a scalar part and a non-scalar part,
denoted, respectively, by
$$x_0={\rm {Sc}}(x), \quad \sum_{k=1}^{n-1}[x]_k={\rm{NSc}}(x). $$

We define the norm of $x\in {\mathcal{CL}}_{0, n}$ to be $|x|=\sqrt{{\rm Sc}(x\bar{x})}=\sqrt{{\rm Sc}(\bar{x}x)}=\sqrt{\sum_{k=0}^{n}\sum_{|A|=k} |x_A|^2}$.

Define the generalized Cauchy-Riemann operator by
$$D=\frac{1}{2}(\frac{\partial}{\partial x_0}+\sum_{k=1}^{n}\frac{\partial}{\partial x_k}{\e}_k)$$
and the related conjugate operator by
$$\bar{D}=\frac{1}{2}(\frac{\partial}{\partial x_0}-\sum_{k=1}^{n}\frac{\partial}{\partial x_k}{\e}_k).$$
Clearly, the Laplace operator $\Delta$ in the $(n+1)$-D space satisfies $\Delta_x=4D\bar{D}=4\bar{D}D$.

Any function $g: {\R}^{n+1}\rightarrow {\mathcal {CL}}_{0, {n}}$ is called  a Clifford-valued function.

\begin{definition}
$g$ is said to be left-monogenic in domain $\Omega\subseteq {\R}^{n+1}$, if $g\in C^1(\Omega)$ and $Dg=0$. If $g$ is left-monogenic, then we call $\bar{D}g$ the left-derivative of $g$.
\end{definition}

\begin{remark} For Clifford monogenic functions we have the following simple conclusions. 1. If $g$ is left-monogenic, then $\displaystyle{\bar D}g=\frac{\partial}{\partial x_0}g;$ and
2. If $g$ is left-monogenic, then $g$ and all its components are harmonic.
\end{remark}

In the Clifford analysis setting, we prove the identical relation between the gradient form and the \emph{monogenic derivative} form quantity.

\begin{theorem}
Let $u(\ux, y)\in {\mathcal D}(\R_+^{n+1}), n\ge 1$, $F$ be the associated canonical left-monogenic function whose scalar part is $u,$ and the non-scalar part has zero value at $\ux=0.$ Then we have
\begin{eqnarray}\label{eqy3}
\|u\|_{\D(\R_+^{n+1}),*}^2&:=&\int_{\R_+^{n+1}}|\nabla u(\ux,y)|^2d\ux dy\nonumber\\
&=&\int_{\R_+^{n+1}}|{\bar D}F(\ux, y)|^2d\ux dy.
\end{eqnarray}
\end{theorem}

To prove the theorem we need, as preparation, the following two lemmas.

\begin{lemma}$^{\mbox{\scriptsize \cite{DMQ}}}$

Let $x_0=y>0$ and $F(\ux, y)=u(\ux, y)+\sum_{k=1}^n v_k(\ux, y)\e_k$ be left-monogenic in $\R_+^{n+1}$, then we have the Cauchy-Riemann equations
\begin{eqnarray}\label{eqy2.4}
DF(\ux, y)=0 \Longleftrightarrow \left\{
\begin{array}{lll}
\vspace{0.2cm}
\displaystyle\frac{\partial u}{\partial y}=\sum_{k=1}^n\frac{\partial v_k}{\partial x_k} \\
\vspace{0.2cm}
\displaystyle\frac{\partial u}{\partial x_k}=-\frac{\partial v_k}{\partial y}, k=1, \cdots, n \\
\vspace{0.2cm}
\displaystyle\frac{\partial v_k}{\partial x_j}=\frac{\partial v_j}{\partial x_k}, j \not= k.
\end{array}
\right.
\end{eqnarray}
\end{lemma}

\begin{lemma}$^{\mbox{\scriptsize \cite{S,DMQ}}}$\label{Cauchy}
Let $u(\ux, y)\in h^2(\R_+^{n+1})$ with the non-tangential boundary value $f(\ux)$, then $F(\ux, y)=u(\ux, y)+\sum_{k=1}^n v_k(\ux, y)\e_k$ is left-monogenic in $\R_+^{n+1}$,
where $v_k(\ux, y)=R_k(u)$ is the $k$-th Rieze transformation of $f$ with the expressions
\begin{eqnarray}\label{eqyy5.8}
v_k(\ux, y)&=&\int_{\R^n}\frac{it_k}{|\ut|}e^{-2\pi {\bf i} \ut\cdot\ux}e^{-2\pi|\ut|y}\hat{f}(\ut)d\ut\nonumber\\
&=&\frac{2}{\omega}\int_{\R^n}\frac{(t_k-x_k)}{(y^2+|\ut-\ux|^2)^{\frac{n+1}{2}}}{f}(\ut)d\ut\nonumber.
\end{eqnarray}
Moreover, $$F(\ux, y)=u(\ux, y)+\sum_{k=1}^n v_k(\ux, y){\e_k}=\frac{1}{\omega_n}E(\cdot, y)\ast f(\ux),$$
where $$E(\ut, y)=\frac{\overline{y+\ut}}{|y+\ut|^{n+1}}$$
is the Clifford Cauchy kernel being monogenic from both the left and right sides.
\end{lemma}

\bigskip
{\bf Proof of (\ref{eqy3}): }If $F(\ux, y)=u(\ux, y)+\sum_{k=1}^n v_k(\ux, y)\e_k$ is left-monogenic in $\R_+^{n+1}$, using (\ref{eqy2.4}) we have

\begin{eqnarray*}
{\bar D}F(\ux, y)&=&\frac{\partial}{\partial y}F(\ux, y)\\
&=&\frac{\partial}{\partial y}[u(\ux, y)+\sum_{k=1}^n v_k(\ux, y)\e_k]\\
&=&\frac{\partial u}{\partial y}+\sum_{k=1}^n \frac{\partial v_k}{\partial y}\e_k\\
&=&\frac{\partial u}{\partial y}-\sum_{k=1}^n \frac{\partial u}{\partial x_k}\e_k.
\end{eqnarray*}

Thus we have
\begin{eqnarray*}
\int_{\R_+^{n+1}}|\nabla u(\ux,y)|^2d\ux dy&=&\int_{\R_+^{n+1}}[(\frac{\partial u}{\partial y})^2+\sum_{k=1}^{n}(\frac{\partial u}{\partial x_k})^2]d\ux dy\\
&=&\int_{\R_+^{n+1}}|{\bar D}F(\ux, y)|^2d\ux dy.
\end{eqnarray*}
This completes the proof of (\ref{eqy3}).

Ahlfors' involvement is to use the quantity defined in one complex variable
\begin{eqnarray}\label{AQ} \frac{1}{2}\int_{\partial \D}\overline{F}F^\prime \frac{dz}{\bf i}\end{eqnarray} together with the residue computation
  in proving the identical relation between the double integral quantity with the other quantities of the Douglas' result (\cite{Ahlfors}). The Douglas' original proof uses the Fourier series expansion of the gradient of the harmonic function (\cite{Douglas}). Validity of the proofs lays on the identical relation between the holomorphic derivative quantity and the Ahlfors' quantity (\ref{more}). In \cite{YQ} by using Stokes' formula the identity (\ref{more}) is generalized to domains in $\R^{n+1}:$
\begin{eqnarray}\label{AH} \int_\Omega |\overline{D}F|^2dV=\frac{1}{2}{\rm Sc}\int_{\partial \Omega}\overline{F}d\sigma \overline{D}F.\end{eqnarray}
The work \cite{YQ} devotes to the balls case in which complicated computation on the left-hand-side of (\ref{AH}) involving spherical harmonic decompositions leads to the double integral form quantity.  With the $\R_+^{n+1}$ case, none of the spherical harmonic analysis methods, nor the complex residue computation that Ahlfors used are available.  We eventually prove the identical relation between the double integral quantity and the Ahlfors form quantities through computation of the right-hand-side of (\ref{AH}). Our treatment using Clifford analysis  can be summarized as computing the integral of a specific left monogenic Hardy $H^1$ function on the upper-half space over the boundary $\R^{n}$ and achieving the desired zero value through using the $\lq\lq$lifting up" method.
In all, this adds the Ahlfors form quantity identical to $\|u\|_{\D(\R_+^{n+1}),*}^2$ as in (\ref{addeqy1}). We have
\begin{theorem}\label{DS}
Let $u(\ux, y)\in {\mathcal D}(\R_+^{n+1})$ with the boundary function $f(\ux)$, then we have
\begin{eqnarray}\label{eqy4}
\|U\|_{\D(\R_+^{n+1}),*}^2&:=&\int_{\R_+^{n+1}}|\nabla u(\ux,y)|^2d\ux dy\nonumber\\
&=&\frac{1}{\omega_n}\int_{\R^n}\int_{\R^n}\frac{|f(\underline{t_1})-f(\underline{t_2})|^2}{|\underline{t_1}-\underline{t_2}|^{n+1}}d\underline{t_1}d\underline{t_2}\nonumber.
\end{eqnarray}
\end{theorem}
{\bf Proof: } For $y>0$, recalling Lemma \ref{Cauchy}, we have
\begin{eqnarray}\label{eqyy5.9}
F(\ux, y)&=&u(\ux, y)+\sum_{k=1}^{n}v_k(\ux, y)\e_k\nonumber\\
&=&\frac{2}{\omega_n}\int_{\R^n}\frac{y+\ut-\ux}{(y^2+|\ut-\ux|^2)^{\frac{n+1}{2}}}f(\ut)d\ut\nonumber\\
&=&\frac{2}{\omega_n}\int_{\R^n}\frac{\overline{y+\ux-\ut}}{(y^2+|\ut-\ux|^2)^{\frac{n+1}{2}}}f(\ut)d\ut.
\end{eqnarray}
Due to the Cauchy integral expression $F(\ux, y)$ is left-monogenic in $\R_+^{n+1}$.
Using (\ref{AH}) and Stokes' formula (cf., \cite{YQ}), we have
\begin{eqnarray*}
&&\int_{\R_+^{n+1}}|\nabla u(\ux,y)|^2d\ux dy\\
&=&\int_{\R_+^{n+1}}|{\bar D}F(\ux, y)|^2d\ux dy\\
&=&-\frac{1}{2}\lim_{y\rightarrow 0^+}\int_{\R^{n}}\bar{F}(\ux, y)d\ux \bar{D}F(\ux, y)\\
&=&-\frac{1}{2}\lim_{y\rightarrow 0^+}\int_{\R^{n}}\bar{F}(\ux, y)d\ux \frac{\partial}{\partial y}F(\ux, y).
\end{eqnarray*}
Taking partial derivative with respect to $y$ to the left-monogenic function $F,$ we get
$$\frac{\partial }{\partial y}F(\ux, y)=\frac{2}{\omega_n}\int_{\R^n}\frac{|\ut-\ux|^2-ny^2-(n+1)y(\ut-\ux)}
{(y^2+|\ut-\ux|^2)^{\frac{n+3}{2}}}f(\ut)d\ut.$$
The change of the orders of taking the derivative and the integration is due to the Lebesgue Dominated Convergence Theorem. The function $\frac{\partial }{\partial y}F(\ux, y)$ is still left-monogenic but with one more order tending to zero as $y$ tending to $\infty.$

For fixed $y$ and $\ut$ the partial derivative of the Cauchy kernel, viz.,
\[\frac{|\underline{t}-\underline{\cdot}|^2-ny^2-(n+1)y(\ut-\underline{\cdot})}
{(y^2+|\ut-\underline{\cdot}|^2)^{\frac{n+3}{2}}},\]
can be viewed as the boundary value of a Hardy $H^1$ function in $\R_+^{n+1} $. Therefore, the Cauchy integra formula may  be used to get
\begin{eqnarray}\label{eqyy34}
&&\int_{\R_+^{n+1}}|\nabla u(\ux,y)|^2d\ux dy=-\frac{1}{2}\lim_{y\rightarrow 0^+}\int_{\R^{n}}\bar{F}(\ux, y)d\ux \frac{\partial}{\partial y}F(\ux, y)\nonumber\\
&=&-\frac{1}{2}\lim_{y\rightarrow 0^+}\frac{2}{\omega_n}\int_{\R^{n}}\int_{\R^n}\frac{y+\ux-\underline{t_2}}
{(y^2+|\underline{t_2}-\ux|^2)^{\frac{n+1}{2}}}f(\underline{t_2})d\underline{t_2}
d\ux\nonumber\\
&\times &\frac{2}{\omega_n}\int_{\R^n}\frac{|\ux-\underline{t_1}|^2-ny^2+(n+1)y(\ux-\underline{t_1})}
{(y^2+|\underline{x}-\underline{t_1}|^2)^{\frac{n+3}{2}}} f(\underline{t_1})d\underline{t_1}\nonumber\\
&=&-\frac{2}{\omega_n}\lim_{y\rightarrow 0^+}\int_{\R^{n}}\int_{\R^{n}}f(\underline{t_2})f(\underline{t_1})
d\underline{t_1}d\underline{t_2}\nonumber\\
&\times &\frac{1}{\omega_n}\int_{\R^n}E((\underline{t_2}+y)-\ux )
d\ux\frac{|\ux-\underline{t_1}|^2-ny^2+(n+1)y(\ux-\underline{t_1})}
{(y^2+|\ux-\underline{t_1})|^2)^{\frac{n+3}{2}}} \nonumber\\
&=&-\frac{2}{\omega_n}\lim_{y\rightarrow 0^+}\int_{\R^{n}}\int_{\R^{n}}f(\underline{t_1})f(\underline{t_2})
\frac{|y+(\underline{t_2}-\underline{t_1})
|^2+y^2+(n+1)y(\ut_2-\underline{t_1})}
{(y^2+|y+(\underline{t_2}-\underline{t_1})|^2)^{\frac{n+3}{2}}}
d\underline{t_1}d\underline{t_2},\nonumber
\end{eqnarray}
where the last step is by replacing $\ux$ with $y+\underline{t_2}$ through using the Cauchy integral formula.
Now we show that for any fixed $y>0$ and $\underline{t_1}\in \R^{n},$
\begin{eqnarray}\label{eqyy35} \int_{\R^{n}}
\frac{|y+(\underline{t_2}-\underline{t_1})
|^2+y^2+(n+1)y(\underline{t_2}-\underline{t_1})}
{(y^2+|y+(\underline{t_2}-\underline{t_1})|^2)^{\frac{n+3}{2}}}d\underline{t_2}
=0.\end{eqnarray}

We note that the integrand is a left-monogenic function, and, in fact, belongs to the Hardy $H^1(\R_+^{n+1}).$  Therefore, the integrals over all the hyper-planes parallelling with the boundary Euclidean plane $\R^n$ are numerically equal to each other. As a consequence, the integral region for the variable $t_2$ can be lifted up to become integral over the hyper-plane $h+t_2, h>0.$ That is,  for any $h>0,$
\begin{eqnarray*} & &\int_{\R^{n}}\frac{|y+(\underline{t_2}-\underline{t_1})
|^2+y^2+(n+1)y(\underline{t_2}-\underline{t_1})}
{(y^2+|y+\underline{t_2}-\underline{t_1}|^2)^{\frac{n+3}{2}}}d\underline{t_2}\\
&=&\int_{\R^{n}}
\frac{|y+h+(\underline{t_2}-\underline{t_1})
|^2+y^2+(n+1)y[h+(\underline{t_2}-\underline{t_1})]}
{(y^2+|y+h+\underline{t_2}-\underline{t_1}|^2)^{\frac{n+3}{2}}}d(\underline{t_2}-\ut_1)\\
&=&\int_{\R^{n}}
\frac{|y+h+\ut|^2+y^2+(n+1)y(h+\ut)}
{(y^2+|y+h+\ut|^2)^{\frac{n+3}{2}}}d\underline{t}. \end{eqnarray*}
Letting $h\to +\infty,$ we obtain the desired relation (\ref{eqyy35}).

As a consequence we have
\begin{eqnarray}\label{zero} \int_{\R^{n}}f^2(\ut_1)\int_{\R^{n}}
\frac{|y+(\underline{t_2}-\underline{t_1})
|^2-ny^2+(n+1)y[y+(\underline{t_2}-\underline{t_1})]}
{(y^2+|y+\underline{t_2}-\underline{t_1}|^2)^{\frac{n+3}{2}}}d\underline{t_2}d\underline{t_1}
=0.\end{eqnarray}

Owing to the symmetry properties, (\ref{zero}) also holds when $f^2(\underline{t_1})$ is replaced by $f^2(\underline{t_2}).$ We finally have

\begin{eqnarray}
& &-\frac{2}{\omega_n}\lim_{y\rightarrow 0^+}\int_{\R^{n}}\int_{\R^{n}}f(\underline{t_1})f(\underline{t_2})
\frac{|y+(\underline{t_2}-\underline{t_1})
|^2-ny^2+(n+1)y[y+(\underline{t_2}-\underline{t_1})]}
{(y^2+|y+(\underline{t_2}-\underline{t_1})|^2)^{\frac{n+3}{2}}}
d\underline{t_1}d\underline{t_2}.\nonumber\\
&=&\lim_{y\rightarrow 0^+}\frac{1}{\omega_n}\int_{\R^{n}}\int_{\R^{n}}
\frac{|f(\underline{t_1})-f(\underline{t_2})|^2\left(|y+(\underline{t_2}-\underline{t_1})
|^2-ny^2+(n+1)y[y+(\underline{t_2}-\underline{t_1})]\right)}
{(y^2+|y+(\underline{t_2}-\underline{t_1})|^2)^{\frac{n+3}{2}}}
d\underline{t_1}d\underline{t_2} \nonumber\\
&=&\lim_{y\rightarrow 0^+}\frac{1}{\omega_n}\int_{\R^{n}}\int_{\R^{n}}
\frac{|f(\underline{t_1})-f(\underline{t_2})|^2}
{(2y^2+|\underline{t_2}-\underline{t_1}|^2)^{\frac{n+1}{2}}}
d\underline{t_1}d\underline{t_2}\nonumber\\
& &\qquad +\lim_{y\rightarrow 0^+}\frac{1}{\omega_n}\int_{\R^{n}}\int_{\R^{n}}
\frac{|f(\underline{t_1})-f(\underline{t_2})|^2
(n+1)y(\underline{t_2}-\underline{t_1})}
{(2y^2+|\underline{t_2}-\underline{t_1}|^2)^{\frac{n+3}{2}}}
d\underline{t_1}d\underline{t_2}. \nonumber
\end{eqnarray}

The second integral is equal to zero for any $y>0$ due to the anti-symmetry property of the integrand. Owing to the Levi Theorem, the first limiting integral is equal to

\begin{eqnarray}
\frac{1}{\omega_n}\int_{\R^{n}}
\int_{\R^{n}}\frac{|f(\underline{t_1})-f(\underline{t_2})|^2}
{|\underline{t_1}-\underline{t_2}|^{n+1}}d\underline{t_1}d\underline{t_2}. \nonumber
\end{eqnarray}

This completes the proof.

In all, we have

\begin{theorem}\label{lastth}
Let $u(\ux, y)\in {\mathcal D}(\R_+^{n+1})$ with the boundary function $f(\ux)$ and $F(\ux, y)$ be the canonical left-monogenic function in $\R_+^{n+1}$ with the scalar part $u(\ux, y)$, then we have
\begin{eqnarray}\label{eqy4}
\|u\|_{\D(\R_+^{n+1}),*}^2&:=&\int_{\R_+^{n+1}}|\nabla u(\ux,y)|^2d\ux dy\nonumber\\
&=&2\pi\int_{\R^{n}}|\ut||\hat{f}(\ut)|^2 d\ut\nonumber\\
&=&\int_{\R_+^{n+1}}|{\bar D}F(\ux, y)|^2d\ux dy\nonumber\\
&=&\frac{1}{\omega_n}\int_{\R^n}\int_{\R^n}\frac{|f(\underline{t_1})-f(\underline{t_2})|^2}{|\underline{t_1}-\underline{t_2}|^{n+1}}d\underline{t_1}d\underline{t_2}.\nonumber
\end{eqnarray}
\end{theorem}

In particular, when $n=1$, we re-obtain Theorem \ref{eqyy1}.

\begin{corollary}\label{last}
Let $u(z)\in {\mathcal D}({\mathbb C}_+)$ with the boundary function $f(x)$ and $F(z)$ be the associated holomorphic function in ${\mathbb C_+}$ with the real part $u(z)$, then we have
\begin{eqnarray}\label{eqy4}
\|u\|_{\D({\mathbb C}_+),*}^2&:=&\int_{{\mathbb C}_+}|\nabla u(z)|^2dx dy\nonumber\\
&=&\int_{{\mathbb C}_+}|F'(z)|^2dx dy\nonumber\\
&=&2\pi\int_{-\infty}^{+\infty}|t||\hat{f}(t)|^2 dt\nonumber\\
&=&\frac{1}{2\pi}\int_{-\infty}^{+\infty}
\int_{-\infty}^{+\infty}\frac{|f(x_1)-f(x_2)|^2}{|x_1-x_2|^2}dx_1dx_2\nonumber.
\end{eqnarray}
\end{corollary}

\bigskip
\subsection{The quaternionic upper-space case}

\
\bigskip

Using the same methods as in the above subsection, we can obtain the analogous results in the quaternionic upper-space. Denote by $\H_+$ the upper-half space of the non-commutative divisible field of quaternionic numbers. With the previously defined and self-explanatory notations we have, in writing  $\H_+=\R^{3}\times (0,\infty).$
\begin{theorem}\label{lastth}
Let $u(\uq, y)\in {\mathcal D}(\H_+)$ with the boundary function $f(\uq)$ and $F(\uq, y)$ the canonical left-regular function in $\H_+$ with the real part $u(\uq, y)$, then we have
\begin{eqnarray}\label{eqy4}
\int_{\H_+}|\nabla u(\uq,y)|^2d\uq dy\nonumber&=&\int_{\H_+}|{\bar D}F(\uq, y)|^2d\uq dy\nonumber\\
&=&2\pi\int_{\R^{3}}|\ut||\hat{f}(\ut)|^2 d\ut\nonumber\\
&=&\frac{1}{2\pi^2}\int_{\R^3}\int_{\R^3}\frac{|f(\underline{t_1})-f(\underline{t_2})|^2}{|\underline{t_1}-\underline{t_2}|^{4}}d\underline{t_1}d\underline{t_2}.\nonumber
\end{eqnarray}
\end{theorem}

\subsection{Balls cases revised}

\
\bigskip

For the Quaternionic unit ball case, as for the unit disc, we have the whole set of identical relations as in Theorem \ref{WZ} and Theorem \ref{lastth}. With a little abuse of notation, denote by $B_4$ the quaternionic unit ball and $S^3$ the quaternionic unit sphere. There holds
\begin{theorem}$^{\mbox{\scriptsize \cite{YQ}}}$
Let $F$ be left-regular in $B_4$ and $u$ its real part with square-integrable boundary value $f$ in $S^3.$ Then there holds
\begin{eqnarray*}
\int_{B_4}|\nabla u|^2dV&=&\int_{B_4}|\bar{D}F|^2 dV\\
&=&\frac{1}{2\pi^2}\int_{S^{3}}\int_{S^{3}}\frac{|f(q_1)-f(q_2)|^2}{|q_1-q_2|^4}
 dS_{\underline{\eta_1}}dS_{\underline{\eta_2}}\nonumber\\
&=&\sum_{k=1}^{+\infty}k\sum_{j=1}^{a^4_k}|b_j|^2,
\end{eqnarray*}
where $b_j$'s are the coefficients of the Fourier-Laplace series of $f$ which are given by formula (\ref{eqY14}) for $n=4.$
\end{theorem}

In summary, for all the cases consisting of upper-half Euclidean space of any dimension, the unit disc case in the complex plane, the quaternionic four dimensional upper-half space and the quaternionic four dimensional ball, there exists, at least morally (with, perhaps, the Fourier series is replaced by Fourier transformation), the same set of identical energy quantities, namely the gradient, the holomorphic or monogenic derivative, the Fourier series or integration, and the double singular integral, as well as the Ahlfors' forms to describe the semi-Dirichlet norm of each context. The only exceptional cases are the unit balls $B_{n+1}, n>1.$ We recall that in
 (\cite{YQ}) we show for $n>1,$
\[\int_{B_{n+1}}|\nabla U(\ux)|^2d\ux\not=\int_{B_{n+1}}|\bar{D} F(\ux)|^2d\ux.\]
Precisely, with the notations that have been using, but with replacement of $n$ by $n+1,$  and consequently, $n>2$ by $n>1$ in Theorem 4.8 of \cite{YQ} (the latter treating only $n>2$ while $n=2$ being a different formula, see Remark 4.1), there holds
\begin{eqnarray*}
\int_{B_{n+1}}|\nabla u(\ux)|^2d\ux&=&
\sum_{k=1}^\infty k\sum_{j=1}^{a^{n+1}_k}|b_j|^2\\
&=&\frac{1}{\omega_n}\int_{S^n}\int_{S^n}\frac{|f(\underline{t_1})-f(\underline{t_2})|^2}
{|\underline{t_1}-\underline{t_2}|^{n+1}}d\underline{t_1}d\underline{t_2},
\end{eqnarray*}
and, in contrast,
\begin{eqnarray*}
\int_{B_{n+1}}|\bar{D} F(\ux)|^2d\ux&=&\frac{1}{2}{\rm Sc}\int_{S^n}\overline{F}d\sigma \bar{D} F\\
&=&\frac{1}{2}\sum_{k=1}(k+\frac{n}{2})\sum_{j=1}^{a^{n+1}_k}|b_j|^2.
\end{eqnarray*}

Although in the case the Ahlfors quantity is not identical with the Douglas quantities,
there is, however, an equivalence relation:
\begin{eqnarray*}
\int_{B_{n+1}}|\bar{D} F(\ux)|^2d\ux&=&
\frac{1}{2}\int_{B_{n+1}}|\nabla u(\ux)|^2d\ux+\frac{n}{4}(\|u\|_{h^2}^2-|u(0)|^2)\\
&\leq&\frac{1}{2}\|u\|_{\D(B_n),*}^2+\frac{n}{4}\|u\|_{\D(B_n)}^2\\
&\backsimeq &\|u\|_{\D(B_n)}^2.
\end{eqnarray*}
So, in the balls case merely the Ahlfors quantity is sufficient to describe the norm of the Dirichlet space.\\

We note that although the Clifford algebra ${\mathcal {CL}}_{0, {2}}$ is identical with the quanternic field $\H$ they are different regarding the Douglas-Ahlfors identities: the former treats $\R^3$ and the latter treats, in a sense, $\R^4.$

\centerline{Acknowledgement}
The authors wish to sincerely thank Prof Wei, Hua-ying, for her kind and inspiriting  communication in relation to this study.

\end{document}